\documentclass[a4paper,11pt]{amsart}
\usepackage[T1]{fontenc}
\usepackage[utf8]{inputenc}
\usepackage[foot]{amsaddr}

\usepackage[pdftex]{hyperref}
\usepackage[margin=2.5cm]{geometry}
\usepackage[english]{babel}
\usepackage{amsmath,amssymb,verbatim,mathrsfs,latexsym,paralist}
\usepackage[dvipsnames]{xcolor}
\usepackage{graphicx}
\usepackage{subcaption}
\usepackage{epstopdf}
\usepackage{indentfirst}
 \usepackage{amsthm}
 \allowdisplaybreaks[4]
 \usepackage{tikz}
\usepackage{caption}
\usepackage[T1]{fontenc}
\usepackage[all]{xy}
\usepackage{xcolor}
\usepackage{soul}

\newcommand{\Z}{\mathbb{Z}}
\newcommand{\C}{\mathbb{C}}

\newcommand{\mrm}{\mathrm}

\newcommand{\rk}{\mathrm{rk}}
\newcommand{\Ext}{\mrm{Ext}}
\newcommand{\Hom}{\mrm{Hom}}
\newcommand{\Span}{\mrm{Span}}

\numberwithin{equation}{section}
\newtheorem{theorem}{Theorem}[section]
\newtheorem*{theorem*}{Theorem}

\newtheorem{proposition}[theorem]{Proposition}

\newtheorem{lemma}[theorem]{Lemma}
\newtheorem{corollary}[theorem]{Corollary}
 
\newtheorem{defn}[theorem]{Definition}  
\newtheorem{example}[theorem]{Example}   

\newtheorem*{lemma*}{Lemma}
\newtheorem*{defn*}{Definition} 
\newtheorem{question}[theorem]{Question}

\theoremstyle{remark}
\newtheorem{rmk}[theorem]{Remark}

\title{Regularity in Linear Degenerations of Flag Variety}
\author[Sabino Di Trani]{Sabino Di Trani}
\address{ Dipartimento di Matematica ``Guido Castelnuovo'', Sapienza - Universit\`a di Roma.
}
\email{sabino.ditrani@uniroma1.it }
\address{\emph{The author has been partially supported by GNSAGA - INDAM group.}}
\address{ORCID id: https://orcid.org/0000-0002-6651-558X}
\begin{document}

\maketitle
\textbf{Abstract:} In this article we investigate the regularity properties of linear degenerations of flag varieties. We classify the linear degenerations of (partial) flag varieties that are smooth. Furthermore, we study the singular locus of irreducible degenerations and provide estimates for its dimension. We also introduce a new stratification of the total space of representations. Within each stratum, we identify the loci corresponding to flat and flat irreducible degenerations. As a consequence of our results, we show that irreducible linear degenerations are normal varieties.
\section{Introduction}
Since their introduction in the 19th century, flag varieties and their Schubert varieties have been central in a great number of problems in combinatorics, topology, algebraic geometry and Lie theory. Research in this area has expanded significantly over the past century, yielding many important results across various mathematical fields.
In the last 20 years, a generalization of these objects—the linear degenerations of flag variety—has been introduced in \cite{CIFFFR20} and \cite{CIFFFR17}. These varieties constitute a broader class that, in many respects, preserves the desirable properties of classical flag varieties.

Given two positive integers $m > n$ and a vector of natural numbers $\mathbf{d}=(d_1, \dots, d_n)$ such that $d_i < d_{i+1}$ and $d_n<m$, consider a family $\mathbf{f}=(f_1, \dots, f_{n-1})$ of linear endomorphisms of~$\C^{m}$.
\begin{defn}A \textbf{linear degeneration} ${\mathcal{F}l}^\mathbf{f}_\mathbf{d}(\C^{m})$ of the partial flag variety ${\mathcal{F}l}_\mathbf{d}(\C^{m})$ associated to the family $\mathbf{f}$ is the algebraic variety of $n$-tuples $(V_1, \dots , V_n)$ of subspaces of $\C^m$ such that $\dim V_i = d_i$ and $f_i V_i\subset V_{i+1}$. 
\end{defn}
Motivated by their multiple connections with the representation theory of algebraic groups, Feigin initially investigated a specific case of degenerations of the complete flag variety ${\mathcal{F}l}(\C^{n+1})$, where the rank of each map $f_i$ is equal to $n$ and the kernels of these maps are linearly independent. These are known as degenerate flag varieties or Feigin degenerations (see \cite{Feigin11}, \cite{Feigin12}). 
In subsequent years, the geometric properties of these objects have been extensively studied using quiver representation theory and combinatorics (see \cite{CIFFFR20}, \cite{CIFFFR17}, \cite{CIFR12}, \cite{CIFR13}, \cite{CIFR17}). Furthermore, in \cite{CIFFFR20} and \cite{CIFFFR17} most of the known results for degenerate flags have been generalized to degenerations obtained considering any $n-1$-tuple of endomorphisms.

In \cite{CIL} it is proved that Feigin degenerations can be realized as Schubert varieties in suitable parabolic flag varieties. This result has been extended to the class of PBW linear degenerations in \cite{CIFFFR20} and \cite{CIFFFR17}. Furthermore, a characteristic free proof that degenerate flags can be relized as (parabolic) Schubert varieties has been provided in~\cite{CILL}.

A fundamental question in the study of Schubert varieties concerns their smoothness and the description of their singular loci. 
Starting from the celebrated pattern avoidance theorem, proved by Lakshmibai and Sandhya in~\cite{LS}, a significant number of works have been devoted to this topic (see \cite{SS} for a comprehensive reference). 

This motivated us to study the smooth and singular loci in linear degenerations, aiming to determine whether they present desirable properties analogous to those found in Schubert varieties.

Surprisingly, only a few results in this direction are already known. 
A first result has been proved in \cite{CIFR13}, where it is shown that  Feigin degenerations are regular in codimension 2, i.e. the singular locus is a closed subvariety of codimension at least 3. 
Furthermore, using the GKM-theory techniques introduced in~\cite{MA}, the smooth locus of PBW linear degenerations has been described in \cite{SDT} as a union of orbits for the action of a suitable Borel group.

Beyond these special cases, the regularity properties of linear degenerations remain largely unexplored. This paper aims to contribute, albeit partially, to fill this gap in the literature.

\subsection{A Geometric Point of View}
Our study adopts a geometric perspective, focusing on the fact that for a fixed dimension vector $\mathbf{d} = (d_1< \dots< d_n)$, linear degenerations of the partial flag variety ${\mathcal{F}l}_\mathbf{d}(\C^{m})$ can be described as fibers of a suitable map. 
Consider the two algebraic varieties  
\[
R := \mathrm{Hom}(\C^{m}, \C^{m})^{n-1}, \qquad Z := \mathrm{Gr}_{d_1}(\C^{m}) \times \dots \times \mathrm{Gr}_{d_n}(\C^{m}),
\]  
where $R$ is endowed with its natural affine space structure. The set  
\[
Y = \{ (\mathbf{f}, \mathbf{V}) \, | \; \mathbf{f} = (f_1, \dots, f_{n-1}) \in R, \, \mathbf{V} = (V_1, \dots, V_n) \in Z, \, f_i V_i \subset V_{i+1} \}
\]  
is a smooth irreducible closed subvariety of $R \times Z$ and is called in~\cite{CIFFFR20} and \cite{CIFFFR17} the \emph{Universal Linear Degeneration of the Flag Variety} ~${\mathcal{F}l}_{\mathbf{d}}(\C^{m})$.
The group $G := \left( GL_m (\C) \right) ^n$ acts on $R$ by change of basis according to the formula  
 \begin{equation}\label{eq:action}  
 (g_1, \dots, g_n) \cdot (f_1, \dots, f_{n-1}) := (g_2 f_1 g_1^{-1}, \dots, g_n f_{n-1} g_{n-1}^{-1}).  
 \end{equation}  
 and on $Z$ componentwise via the classical matrix action on~$\C^{m}$. Consequently, $G$ acts on the product variety $R \times Z$, and the subvariety $Y$ turns out to be a $G$-variety; furthermore the projection $\pi: Y \rightarrow R$ results to be a proper $G$-equivariant morphism for this action.
Many geometric properties of linear degenerations can be deduced by analyzing the properties of $\pi$.
 
As shown in \cite{CIFFFR20} and \cite{CIFFFR17}, there exist two maximal $G$-invariant subsets $$ U_{flat, irr} \subset U_{flat} \subset R,$$ called the flat irreducible locus and flat locus, respectively, such that  
\begin{itemize}  
    \item the restriction of $\pi$ to $\pi^{-1} \left(U_{flat}\right)$ is a flat morphism of algebraic varieties;  
    \item among the $n-1$-tuples in $U_{flat}$, the fiber $\pi^{-1}(\mathbf{f})$ is an irreducible algebraic variety if and only if $\mathbf{f} \in U_{flat, irr}$.  
\end{itemize}  

The locus $U_{flat}$ corresponds to the maximal locus in $R$ such that the morphism $$\pi : \pi^{-1}(U_{flat}) \rightarrow U_{flat}$$ is flat with fibers of minimal dimension.
Degenerations corresponding to elements in $U_{flat}$ are called flat degenerations.
The orbit of $(\mrm{id}, \dots, \mrm{id})$, i.e. the orbit of $n-1$-tuples with fibers isomorphic to the complete flag variety, turns out to be an open dense subset of $U_{flat}$;
since the morphism $\pi$ is flat when restricted to $U_{flat}$, all the degenerations parametrized by this set result to be equidimensional, with dimension equal to $\dim {\mathcal{F}l}_{\mathbf{d}}(\C^{m})$.
\subsection{Main Results:} In this paper, we focus on the following three questions:
\begin{enumerate}
    \item[(Q1)] Can we characterize which degenerations are smooth?
    \item[(Q2)] If $\mathcal{F}l_\mathbf{d}^\mathbf{f} (\C^m)$ is not smooth, is it possible to explicitly describe its singular locus?
        \item[(Q3)] Can conditions be given on the maps $f_1, \dots, f_{n-1}$ such that the singular locus has interesting geometric properties (irreducibility, smoothness, dimension estimates)? 
\end{enumerate}
Our main theorems address the first question by providing a complete characterization of irreducible and smooth linear degenerations. 
\begin{theorem*}
  $\mathcal{F}l_\mathbf{d}^\mathbf{f} (\C^m)$ is an irreducible algebraic variety if and only if for all $i \in \{1, \dots, n-1\}$ either $ d_{i+1} -d_i \geq m-\rk f_i$ or $\rk f_i=0$ . 
\end{theorem*}
\begin{theorem*} A linear degeneration $\mathcal{F}l_\mathbf{d}^\mathbf{f} (\C^m)$ is smooth if and only if for all $i \in \{1, \dots, n-1\}$ either $\rk f_i = 0$ or $\rk f_i = m$.
\end{theorem*}
As a result of the characterization of irreducible linear degenerations, we prove that  
\begin{theorem*} Every irreducible linear degeneration is a normal variety. 
\end{theorem*}
Furthermore, as a consequence of our classification of smooth linear degenerations, we introduce a new stratification of the space $R$. Within each stratum, we identify both a flat and a flat irreducible locus, which we describe in terms of rank sequences. This description generalizes that of $U_{flat}$ and $U_{flat, irr}$ from \cite{CIFFFR20} and \cite{CIFFFR17}.
The union of these flat irreducible loci across all strata forms the complete irreducible locus $U_{irr} \subset R$ of $n-1$-tuples of endomorphisms that parametrize irreducible degenerations.

Subsequently, in order to extend the regularity result proved in \cite{CIFR13} for Feigin degenerations, we focus on flat irreducible degenerations. Using the properties of the map $\pi$, we prove that if $X$ and $X'$ are flat degenerations and $X$ degenerates to $X'$, then  
\[\dim \mrm{Sing}(X') \geq \dim \mrm{Sing}(X),\]
where $\mrm{Sing}(X)$ denotes the singular locus of the variety $X$. 
As a consequence of this argument, we prove that flat irreducible degenerations are regular in codimension 2. By virtue of our classification theorem, we are able to extend this result to all irreducible degenerations:
 \begin{theorem*}
Every irreducible linear degeneration $X$ is regular in codimension 2. Furthermore, if $d_{i+1}-d_i=1$ for all $i \in \{1, \dots, n\}$, the singular locus $\mrm{Sing}(X)$, when non empty, is a closed subvariety of $X$ of codimension 3.
\end{theorem*}
Finally, we pose several open problems (see Questions~\ref{q:1}, \ref{q:2} and \ref{q:3}), concerning the dimension of the singular locus when the conditions on the dimension vector are relaxed.
\subsection{Organization of the Paper:}
Since linear degenerations can be conveniently described as quiver Grassmannians, in Section 2 we recall the main tools coming from the representation theory of quivers. Section 3 is devoted to proving the classification of irreducible and smooth linear degenerations. In Section 4 we investigate the new stratification of the total space of representations and we study flat and flat irreducible loci in these new strata. Finally, Section 5 focuses on the study of the singular locus of certain flat linear degenerations. We first recall some results from algebraic geometry and then apply them to study the regularity of irreducible linear degenerations.

\subsection{Acknowledgments:} Most of the results presented in this article have been obtained during my stay at the Chair of Algebra and Representation Theory at RWTH Aachen University. I would like to thank all members of Ghislain Fourier’s research group for their warm hospitality, and in particular Xin Fang for his availability and support throughout this research. I am also grateful to Giovanni Cerulli Irelli for many insightful discussions related to the topic of this paper. Moreover, I sincerely thank Mattia Talpo and Davide Gori for answering my questions concerning the algebraic geometry tools used in this work. Finally, I would like to thank Martina Lanini for pointing out a mistake in a preliminary version of the proof of Theorem~\ref{thm:singcodim3}.
My stay in Aachen was financially supported by an INdAM Travel Fellowship (Mensilità per mobilità all’estero INdAM).

\section{Quiver Representation Theory Background}
Let $Q$ be a finite quiver with set of vertices $Q_0$ and set of oriented edges $Q_1$. 
A finite dimensional representation $M$ of $Q$ with base field  $\C$ is the datum of a family of finite dimensional complex vector spaces $M=(M_i)_{i \in Q_0}$ and a set of maps $(f_\alpha)_{\alpha \in Q_1 }$ such that, if $\alpha$ is the arrow $i \rightarrow j, \; i,j \in Q_0$, then  $f_\alpha \in \mrm{Hom}_\C(M_{i}, M_{j})$.  For a complete reference on this subject, we refer to \cite{ASS} and \cite{Ki}.
A morphism between two representations $M=(M_i, f_\alpha)$ and $N=(N_i, g_\alpha)$ is a tuple of maps $(\psi_i)_{i \in Q_0}$ such that $ \psi_i$ is a linear map from $M_i$ to $N_i$ and $ \psi_j \circ f_\alpha = g_\alpha \circ \psi_i$ for any edge $\alpha: i \rightarrow j$. Morphisms can be composed componentwise and $\mrm{Rep} Q$ results to be an abelian category, equivalent to the category of finite dimensional left modules over the path algebra of $Q$. 
To each representation $M$ is attached a dimension vector ${\bf{dim M}}=(m_i)_{i \in Q_0}$ such that $m_i=\mrm{dim}_\C M_i$.
A subrepresentation $N$ of $M$ is a sub-object of $M$ in the category $\mrm{Rep}(Q)$. 
\begin{defn}
Let ${\bf{d}}=(d_i)_{i \in Q_0}$ be a dimension vector. The quiver Grassmannian $\mrm{Gr}_{\bf d}(Q, M)$ is the algebraic variety of subrepresentations of $M$ with dimension vector $\bf d$. 
\end{defn}
We are interested in working with representation of equioriented quivers of type $A_n$, i.e. the quiver defined by the data  $Q_0=\{1, \dots, n\}$ and $Q_1=\{(i,i+1) | 1 \leq i < n\}$, we then omit the quiver $Q$ in the quiver Grassmannian notation.
If $Q$ is the equioriented quiver of type $A_n$, the Euler Form is defined by the formula
$$\langle{\bf d},{\bf e}\rangle=\sum_{i=1}^nd_ie_i-\sum_{i=1}^{n-1}d_ie_{i+1}.$$
Since the category $\mrm{Rep}(A_n)$ is hereditary, i.e. $\mrm{Ext}^{\geq 2}(M,N)=0$ for every pair of representations $M,N$, the following formula holds: 
\begin{equation}\label{ExtHom}
\dim_\C \Hom_Q(M,N)-\dim_\C \Ext^1(M,N)=\langle \bf{dim M}, \bf{dim N} \rangle. 
\end{equation} 
By Gabriel's Theorem, the indecomposable representations of the equioriented quiver $A_n$ are all of the form $U_{i,j}$, for $1\leq i\leq j\leq n$, where $U_{i,j}$ is the representation 
$$0\rightarrow\ldots\rightarrow 0\rightarrow \mathbb{C}\stackrel{{\rm id}}{\rightarrow} \ldots \stackrel{{\rm id}}{\rightarrow}\mathbb{C}\rightarrow 0\rightarrow\ldots\rightarrow 0,$$
supported on the vertices with indices between $i$ and $j$.
In particular injective and projective objects in $\mrm{Rep} A_n$ are of the form $I_i=U_{1,i}$ and $P_i=U_{i,n}$, respectively.
With $S_i$ we denote the simple object $U_{i,i}$, supported on the $i$-th vertex of the quiver $A_n$.
For representations of the form $U_{i,j}$ it is easy to compute the dimension of $\Hom$ and $\Ext^1$: 
\small{\[
\dim \Hom_Q(U_{i,j},U_{h,k})=
\begin{cases}
 1 & \text{if} \;h \leq i \leq k \leq j \\
 0 & \text{otherwise.}
\end{cases}
\qquad
\dim \Ext^1(U_{i,j},U_{h,k})=
\begin{cases}
 1 & \text{if} \;i+1 \leq h \leq j+1 \leq k  \\
 0 & \text{otherwise.}
\end{cases}
\]}
\subsection{Linear Degenerations as Quiver Grassmannians}
Let $m$ be a fixed positive integer greater than $n$ and $\mathbf{f}=\{f_1, \dots , f_{n-1}\} $ a $n-1$-tuple of endomorphisms of $\C^{m}$. Consider the representation 
 $$ M_{\bf f} := \quad  \C^m \stackrel{{f_1}}{\longrightarrow} \C^m \stackrel{{f_2}}{\longrightarrow} \ldots  \stackrel{{f_{n-2}}}{\longrightarrow}\C^m \stackrel{{f_{n-1}}}{\longrightarrow}\C^m $$ 
of the equioriented quiver of type $A_n$.
Given a dimension vector $\mathbf{d}=(d_1, \dots, d_n)$, such that $0 < d_i < d_{i+1}< m$ for all $i$, the linear degeneration~${{\mathcal{F}l}^{\bf{f}}_{\bf{d}}}(\C^{m})$ can be realized as the quiver Grassmannian $\mrm{Gr}_{\bf{d}}( M_{\bf{f}})$. 

Linear degenerations of the partial flag variety ${\mathcal{F}l}_{\mathbf{d}}(\C^{m})$ are parametrized by $n-1$-tuples in the space $R:=\oplus_{i=1}^{n-1} \mathrm{Hom}(M_i, M_{i+1}) = \left( \mrm{End} (\C^m) \right)^{n-1}$. The group $G:= \left( GL_m (\C) \right) ^n$ acts on $R$ according to change of basis rule given by formula~\ref{eq:action}.
There is a correspondence between isomorphism classes of representations and $G$-orbit in $R$.
These orbits are uniquely identified by sequences of integers $r_i$ and $ r_{i,j}$ with $ i <j \leq n-1$, where $r_{i} = \rk f_i$ and $r_{i,j}$ is the rank of composition map $f_{j} \circ \dots \circ f_i$. 
The rank sequences also identify the closure ordering on orbits. It we denote the rank sequence identifying the orbit $O_M$ of representation $M$ as $r(M)$, it easy to check that $O_M$ is contained in the closure of  $O_N$ if and only if $r(M) \leq r(N)$ componentwise. In this case we say that $N$ degenerates to $M$ and that the orbit $O_N$ degenerates to $O_M$.
\subsection{Orbit representatives}\label{sec:orbitrep}Let us fix a basis $B=\{v_1, \dots, v_{m}\}$ of $\C^{m} $ and  $J \subseteq \{1, \dots, m\}$. The projection operator $\pi_J: \C^{m} \rightarrow \C^{m}$ (with respect to $B$) is the linear operator 
on $\C^{m} $ 
defined by the rule: 
\[\pi_J(v_i)=\begin{cases} 0, & \text{if } i \in J;  \\ v_i, & \text{ otherwise. } \end{cases}\]
For any family $\mathbf{J}=(J_1, \dots, J_{n-1})$ of subsets of $\{1, \dots, m\}$ we define the degenerate flag variety ${\mathcal{F}l}_{\mathbf{d}}^\mathbf{J}(\C^{m} )$ as the $\mathbf{f}$-degeneration of ${\mathcal{F}l}(\C^{m} )$ obtained considering the family of endomorphisms $\mathbf{f}=(f_1, \dots, f_{n-1})$ such that 
$f_i=\pi_{J_i}$ for every $i$.
Acting with $G$ by change of basis, it is possible to prove that each $G$-orbit has a representative of the form $(\pi_{J_1}, \dots, \pi_{J_{n-1}})$ and consequently each fiber is isomorphic to a linear degeneration ${\mathcal{F}l}_{\mathbf{d}}^\mathbf{J}(\C^{m} )$ for a suitable family of subsets $\bf J$. 
\subsection{Loci in $\mrm{Rep}\label{sec:loci} \,Q$ and Properties of Linear degenerations}
Since $n-1$-tuples in the same orbit correspond to isomorphic linear degenerations, it is natural to investigate whether specific loci in $R$ correspond to degenerations with desirable geometric properties. As an example, the $G$ orbit of $(\mrm{id}, \dots, \mrm{id})$ is the unique dense open orbit in $R$, all the fibers over point of such orbit are isomophic to the partial flag ${\mathcal{F}l}_{\mathbf{d}}(\C^{m})$ and in particular they are irreducible, normal and smooth projective varieties. It is natural to ask wether it is possible to describe all the $\mathbf{f} \in R$ such that  ${{\mathcal{F}l}^{\bf{f}}_{\bf{d}}}(\C^{m})$ is irreducible and, eventually smooth. In the case of linear degenerations of the complete flag variety, in \cite[Theorem 13]{CIFFFR17} it has been shown that a flat linear degeneration ${{\mathcal{F}l}^{\bf{f}}_{\bf{d}}}(\C^{m})$ is a normal variety if and only if $\mathbf{f} \in U_{flat, irr}$. In order to extend these results, it is useful to recall how the loci $U_{flat}$ and $U_{flat, irr}$ can be described, by means of rank sequences.
Consider the rank sequences $\mathbf{r}^1$ and $\mathbf{r}^2$ defined by 
\[r^1_i = m +d_i - d_{i+1}, \qquad r^1_{ij} = m +d_i - d_{j+1},\]
\[r^2_i = m +d_i - d_{i+1} -1, \qquad r^2_{ij} = m +d_i - d_{j+1}-1.\]
Let us denote by $O_{\mathbf{r}^1}$ and $O_{\mathbf{r}^2}$ the orbits identified by $\mathbf{r}^1$ and $\mathbf{r}^2$, respectively.
\begin{theorem}\cite[Theorem A]{CIFFFR20}\label{thm:flatlocus}\begin{enumerate}\item The set $U_{flat}$ is the locus of orbits degenerating to $O_{\mathbf{r}^2}$, \item Let $\mathbf{f}  \in U_{flat}$, the degeneration ${{\mathcal{F}l}^{\bf{f}}_{\bf{d}}}(\C^{m})$ is a reduced locally complete intersection variety,
\item The set $U_{flat,irr}$ is the locus of orbits degenerating to $O_{\mathbf{r}^1}$. \end{enumerate}
\end{theorem}
\subsection{Schubert Quiver Grassmannians}
We recall that the Auslander-Reiten quiver $\Gamma$ of the equioriented quiver of type $A_n$ is described by the following data: the vertices of $\Gamma$ are parametrized by pairs $(i,j)$ with $1 \leq i \leq j \leq n$; its arrows are $(i,j) \rightarrow (i-1, j)$ and $(i,j) \rightarrow (i,j-1)$, whenever the targets are well defined. 

\begin{defn}\cite[Definition 2.1]{CIFR17} A representation $M$ is catenoid if all distinct indecoponsable direct summands  of $M$ belong to an oriented connected path of the Auslander–Reiten quiver.
\end{defn}
Recall that a minimal projective resolution of a representation $M$ is given by
\[ 0 \rightarrow Q \rightarrow P \rightarrow M \rightarrow 0\]
where \[Q=\bigoplus_{1 \leq i \leq n} P_i^{\dim \Hom_Q(M, S_i)}, \qquad P= \bigoplus_{1 \leq i \leq n} P_i^{\dim \Ext^1(M, S_i)} .\]
Such a resolution provides a map of algebraic varieties
\[\mrm{Gr}_\mathbf{d}(M) \rightarrow \mathcal{F}l_{\mathbf{d}+\mathbf{dim Q}}(P)\]
that results to be a closed embedding. In \cite[Theorem 2.2]{CIFR17} it is also proved that the image of this map is $B$-stable if and only if the representation $M$ is catenoid. Under these assumptions, its irreducible components can be identified with certain Schubert varieties, explicitly realizable (see \cite[Section 4]{CIFR17}). A quiver Grassmannian $\mrm{Gr}_\mathbf{d}(M)$ such that $M$ is catenoid is a Schubert quiver Grassmannian.

\section{Irreducibility and Smoothness Condition}\label{sec:smooth}
Consider $\mathbf{f}=(f_1, \dots, f_{n-1}) \in R$ and a dimension vector of positive integers ${\bf d}=(d_1 < \dots < d_n )$ such that $d_n < m$.
In this section, our interest is focused on determining how the irreducibility and smoothness of $Gr_{\mathbf{d}}(M_\mathbf{f}) \simeq \mathcal{F}l_{\mathbf{d}} (\mathbb{C}^m)$ depend on the rank of the maps~$f_i$. 
\begin{rmk}\label{rmk:product} If some of the elements in $\mathbf{f}$ are the null map, $Gr_{\mathbf{d}}(M_\mathbf{f})$ splits as a product of quiver Grassmannians for suitable smaller quivers. More precisely, assume that $f_i=0$ for an $i \in \{1, \dots n-1\}$. Set \[\mathbf{f'}=(f_1, \dots, f_{i-1}, 0, \dots, 0), \qquad \mathbf{d'}=(d_1, \dots, d_{i}, 0, \dots, 0),\]\[ \mathbf{f''}=(0, \dots, 0, f_{i+1}, \dots, f_{n-1}), \qquad \mathbf{d''}=(0, \dots, 0, d_{i+1}, \dots, d_{n}),\]
we have $$M_\mathbf{f} = M_\mathbf{f'} \oplus M_\mathbf{f'}.$$
Consider the two representations \[\widetilde{M_\mathbf{f'}}:= \C^m \stackrel{f_1} \to \C^m \stackrel{f_2} \to \dots \stackrel{f_{i-2}} \to \C^m\stackrel{f_{i-1}} \to \C^m, \qquad    \widetilde{M_\mathbf{f'}}:=\C^m \stackrel{f_{i+1}} \to \C^m \stackrel{f_{i+2}} \to \dots \stackrel{f_{N-2}} \to \C^m\stackrel{f_{n-1}} \to \C^m\] of the equioriented quivers of type $A_i$ and $A_{n-i}$ respectively and the two dimension vectors 
\[\widetilde{\mathbf{d'}}=(d_1, \dots, d_{i}), \qquad \widetilde{\mathbf{d''}}=( d_{i+1}, \dots, d_{n}).\] We have 
\[Gr_{\mathbf{d'}}(M_{\mathbf{f}'}) \simeq Gr_{\widetilde{\mathbf{d'}}}(\widetilde{M_{\mathbf{f}'}}), \qquad Gr_{\mathbf{d''}}(M_{\mathbf{f}''}) \simeq Gr_{\widetilde{\mathbf{d}''}}(\widetilde{ {M_{\mathbf{f}''}} }) \]and consequently $$Gr_{\mathbf{d}}(M_{\mathbf{f}})\simeq Gr_{\widetilde{\mathbf{d'}}}(\widetilde{M_{\mathbf{f}'}}) \times Gr_{\widetilde{\mathbf{d''}}}(\widetilde{ {M_{\mathbf{f}''}} }). $$

\end{rmk}
We now classify irreducible linear degenerations of partial flags. We firstly discuss the case $n=2$. Consider the representation 
\[ M_{f}:= \C^m \stackrel{f} \longrightarrow \C^m.\]
If $f$ is the null map, then $\mrm{Gr}_{\mathbf{d}}(M_f)$ is product of classical Grassmannians and in particular it is an irreducible variety. The other cases are covered by the following Lemma. 
\begin{lemma}\label{lem:basestep} Consider a non zero $f \in End(\C^m)$ and ${\bf{d}}=(d_1, d_2)$, a dimension vector. The quiver Grassmannian $\mrm{Gr}_{\mathbf{d}}(M_f)$ is irreducible if and only if $d_2 - d_1 \geq m - \mrm{rk} f$
\end{lemma}
\proof Since the representation $M$ is catenoid 
the variety $\mrm{Gr}_{\mathbf{d}}(M)$ is a Schubert Quiver Grassmannian. 
Set \[ Q=P_2^{m-\mrm{rk} f}, \qquad P= P_1^m \oplus P_2^{m - \rk f}.\] A minimal projective resolution for $M$ is given by 
\[ 0 \rightarrow Q \rightarrow P \rightarrow M \rightarrow 0.\] 
The dimension vectors $\mathbf{p}=(p_1, p_2)$ and $\mathbf{q}=(q_1, q_2)$ of $P$ and $Q$ respectively, are given by $\mathbf{p}=(m , 2m-\rk f)$ and $\mathbf{q}=(0 , m-\rk f)$. Consequently, the representation $M$ is simple (see~\cite[Definition 4.3]{CIFR17}). Because of \cite[Theorem 4.4]{CIFR17} it is irreducible if and only if $p_1 - d_1 \geq p_2 - d_2$. The latter inequality is equivalent to require that $d_2 - d_1 \geq m - \mrm{rk} f$.
\endproof
To prove the inductive step, we need a preliminary Lemma. 
\begin{lemma}\label{lem:rankcomposition}Let $\mathbf{d}=(d_1 < \dots < d_n)$ be a dimension vector and $\mathbf{f}=(f_1, \dots, f_{n-1})$ be a $n-1$-tuple of endomorphisms. Suppose that $ d_{h+1} - d_h \geq m-\rk f_h$  for every $ h \in \{1, \dots, n-1\}.$
Then $$ r_{hk} = \rk ( f_k \circ \dots \circ f_h) \geq m + d_h - d_{k+1}$$ 
for every $h,k$ such that $h < k$.
\end{lemma}
\proof
To prove the assertion, it is enough to show that  $ r_{hh+1} \geq m + d_h - d_{h+2}$ and then proceed by iteration. The assertion for $k=h+1$ can be obtained by the following computations: 
\begin{align*}
\rk ( f_{i+1} \circ f_i ) &\geq \rk f_i - \dim \left( \mathrm{Ker} f_{i+1} \cap \mathrm{Im} f_i \right) \\
&\geq \rk f_i - (m-\rk f_{i+1}) \\
&= \rk f_i + \rk f_{i+1} -m \\
& \geq m +d_i - d_{i+1} + m + d_{i+1} - d_{i+2} - m \\
& = m +d_i - d_{i+2}.
\end{align*}
\endproof
Recall now that the product of two complex algebraic variety $X_1 \times X_2$ is irreducible if and only if both $X_1$ and $X_2$ are irreducible. 
\begin{theorem}\label{thm:irreducible}
 The quiver Grassmannian $Gr_{\mathbf{d}}(M_\mathbf{f})$ is an irreducible algebraic variety if and only if for all $i \in \{1, \dots, n\}$ either $ d_{i+1} -d_i \geq m-\rk f_i$ or $\rk f_i=0$. 
\end{theorem}
\proof 
Firstly, we prove that if $d_{i+1} -d_i \geq m-\rk f_i$ or $\rk f_i=0$ for all $i \in \{1, \dots, n\}$, then  $Gr_{\mathbf{d}}(M_\mathbf{f})$ is irreducible. Due to Remark~\ref{rmk:product}, we can reduce to prove irreducibility when every $f_i$ is not the null map, i.e. we can suppose that $d_{i+1} -d_i \geq m-\rk f_i$ for all $i$.
As a consequence of Lemma~\ref{lem:rankcomposition} we have that  $ r_{hk}  \geq m + d_h - d_{k+1}$
for every $h,k$ such that $h < k$. According to \cite[Theorem 3]{CIFFFR20}, the quiver Grassmannian $Gr_{\mathbf{d}}(M_\mathbf{f})$ can then be identified with a flat irreducible linear degeneration of a suitable partial flag variety. \\
Let us assume now that $Gr_{\mathbf{d}}(M_\mathbf{f})$ is irreducible.
We proceed by induction on $n$.
If $n=2$ the assertion is estabilished by Lemma~\ref{lem:basestep}.
Suppose now $n > 2$ and consider the representation $M^s$ of the equioriented quiver $A_{n-1}$ defined by the data 
\[M^s_i = \mathbb{C}^m \mbox{ for every } i \in [n-1], \qquad \qquad M^s_{i , i+1} = f_i \mbox{ for every } i \in [n-2].\]
Set $\mathbf{d}^s=(d_1, \dots, d_{n-1})$. Observe that, since $Gr_{\mathbf{d}}(M_\mathbf{f})$ is not empty, then also the quiver Grassmannian $Gr_{\mathbf{d}^s}(M^s)$ is not empty.
The projection map 
\[\begin{matrix}
\pi^s:&Gr_{\mathbf{d}}(M_\mathbf{f}) &{\longrightarrow} &Gr_{\mathbf{d}^s}(M^s)\\
&(V_1, \dots, V_n) &\longrightarrow & (V_1, \dots, V_{n-1})
\end{matrix} \]
is a surjective morphism of algebraic varieties, and consequently if $Gr_{\mathbf{d}}(M_\mathbf{f})$ is irreducible, then $Gr_{\mathbf{d}^s}(M^s)$ also needs to be irreducible. By induction hypothesis, we obtain either $d_i - d_{i+1} \geq m-\rk f_i$ or $f_i = 0$, for all $i \leq n-2$.  Analogously, set \[M^e_i = \mathbb{C}^m \mbox{ for } i \in [n-1] \qquad M^e_{i , i+1} = f_{i+1} \mbox{ for } i \in [n-2], \qquad \mathbf{d}^e=(d_2, \dots, d_n)\]
and consider
\[\begin{matrix}
\pi^e:&Gr_{\mathbf{d}}(M_\mathbf{f}) &{\longrightarrow} &Gr_{\mathbf{d}^e}(M^e)\\
&(V_1, \dots, V_n) &\longrightarrow & (V_2, \dots, V_{n})
\end{matrix} \]
 As in the previous case, $Gr_{\mathbf{d}^e}(M^e)$ need to be irreducible, and by inductive hypothesis we have either $d_{i+1} -d_i \geq m-\rk f_i$ or $f_i = 0$, for all $i$ such that $ 2 \leq i \leq n-1$.
\endproof
\begin{corollary}\label{cor:euler} A linear degeneration $\mathcal{F}l_{\bf{d}} (\mathbb{C}^m)$
is irreducible if and only if it is a product of flat irreducible linear degenerations of partial flag varieties. \end{corollary}
\proof 
If $\mathcal{F}l_{\bf{d}} (\mathbb{C}^m)$  is a product of irreducible flat degenerations of partial flag varieties, then it is an irreducible algebraic variety. Conversely, we can identify $\mathcal{F}l_{\bf{d}} (\mathbb{C}^m)$ with the quiver Grassmannian $Gr_{\mathbf{d}} (M_\mathbf{f})$ and consider the set 
$$I=\{ i \in \{1, \dots n-1\}  \, | \, f_i=0\} \cup \{0\}.$$
If $I=\{i_1 \leq \dots \leq i_k\}$, for $j \in \{ 1, \dots k\}$, set 
\[\mathbf{f}^j=(f_{i_j+1}, \dots, f_{i_{j+1}-1}),\qquad \mathbf{d}^j=(d_{i_j+1}, \dots,d_{i_{j+1}}),\qquad X_j=Gr_{\mathbf{d}^j}(M_{\mathbf{f}^j}).\]
 Because of Remark~\ref{rmk:product}, $$Gr_{\mathbf{d}} (M_\mathbf{f})=X_1 \times \dots \times X_k.$$
As a consequence of Lemma~\ref{lem:rankcomposition}, we have thah $ r_{hk} \geq m + d_h - d_{k+1}$ for every $h,k \in \{i_j+1, \dots, i_{j+1}-1\}$ such that $h< k$. By virtue of \cite[Theorem 3]{CIFFFR20}, for every $j$, the quiver Grassmannian $X_j$ is then isomorphic to a flat irreducible linear degeneration of the partial flag variety $\mathcal{F}l_{\bf{d}^j} (\mathbb{C}^m)$. \endproof

\begin{rmk}\label{rmk:dim}
Recall that the partial flag variety $\mathcal{F}l_{\bf{d}^j} (\mathbb{C}^m)$, has dimension equal to
\[ \langle \mathbf{d}^j, \dim M_{\mathbf{f}^j} - \mathbf{d}^j \rangle = d_{i_j+1}(m-d_{i_j+1})+\sum_{h=i_j+2}^{i_{j+1}} (d_h-d_{h-1})(m -d_h ).\] Since $X_j$ can be identified with a flat irreducible linear degeneration of $\mathcal{F}l_{\bf{d}^j} (\mathbb{C}^m)$ we have that 
\[\dim X_j = \langle \mathbf{d}^j, \dim M_{\mathbf{f}^j} - \mathbf{d}^j \rangle \] \end{rmk}
Recall that, by \cite[Proposition 6]{CR}, if $L $ is a subrepresentation of dimension vector $\mathbf{d}$ of a quiver representation $M$, then the tangent space to  $Gr_\mathbf{d}(M)$ at $L$ is isomorphic as a vector space to $\mrm{Hom}_Q(L, M/L)$. 
As consequence of Euler formula, if $Gr_{\mathbf{d}} (M_\mathbf{f})=X_1 \times \dots \times X_k$ we have that,
\begin{align*} \dim T_L ( Gr_{\mathbf{d}}(M_\mathbf{f}) ) &= \sum_j \dim \mrm{Hom}_Q(L^j, M_{\mathbf{f}^j}/L^j) \\  &=  \sum_{j} \left( \langle \mathbf{d}^j, \dim M_{\mathbf{f}^j} - \mathbf{d}^j \rangle  + \dim \mrm{Ext}^1 (L^j , M_{\mathbf{f}^j}/L^j) \right) \\ &=  \sum_{j} \left( \dim  Gr_{\mathbf{d}^j}(M_{\mathbf{f}^j})  + \dim \mrm{Ext}^1 (L^j , M_{\mathbf{f}^j}/L^j) \right) , \end{align*}
where $L^j$ is the projection of $L$ into $M_{\mathbf{f}^j}$.
According to Remark~\ref{rmk:dim}, $L$ corresponds to a smooth point in the irreducible degeneration  $\mathcal{F}l_\mathbf{d}^\mathbf{f} (\C^m) \simeq Gr_{\mathbf{d}}(M_\mathbf{f})$ if and only if $\mrm{Ext}^1 (L^j , M_{\mathbf{f}^j}/L^j) = 0$ for all $j$. \\
As a consequence of this observation, we classify the linear degenerations of partial flag varieties that are smooth. 
\begin{theorem} \label{thm:smooth} A linear degeneration $\mathcal{F}l_\mathbf{d}^\mathbf{f} (\C^m)$ is smooth if and only if for all $i \in \{1, \dots, n\}$ either $\rk f_i = 0$ or $\rk f_i = m$.
\end{theorem}
\proof If  $\rk f_i = 0$ or $\rk f_i = m$ for every $i$, the degeneration $\mathcal{F}l_\mathbf{d}^\mathbf{f} (\C^m)$ is a product of partial flag varieties, and hence it is irreducible and smooth. To prove the other implication, observe that, since $\mathcal{F}l_\mathbf{d}^\mathbf{f} (\C^m)$ is connected, it is smooth only if it is irreducible. In virtue of Theorem~\ref{thm:irreducible}, we can suppose that either $d_i - d_{i+1} \geq m-\rk f_i$ or $\rk f_i=0$ for all $i \in \{1, \dots, n\}$. Moreover, because of Remark~\ref{rmk:product}, we can reduce to the case of $\rk f_i \neq 0$ for all $i$, i.e. $\mathcal{F}l_\mathbf{d}^\mathbf{f} (\C^m)$ is a flat irreducible linear degeneration.

We prove that if $m > \rk f_i$ for some $i$, then there exists a $n$-tuple of subspaces $F=(F_1, \dots, F_n) \in \mathcal{F}l_\mathbf{d}^\mathbf{f} (\C^m)$ corresponding to a point $L_F$ of $Gr_{\mathbf{d}}(M_\mathbf{f})$ such that  $\mrm{Ext}^1 (L_F , M_{\mathbf{f}}/L_F)  \neq 0$. This point is not smooth in virtue of Remark~\ref{rmk:dim}. Let $h$ be smallest index such that $\rk f_h < m$.  Set $k=m - \rk f_h$. Up to change of basis, it is enough to show our thesis for all the $f_i$ being of the form $\pi_{J_i}$ with respect to a suitable basis $\{v_1, \dots, v_m\}$ of $\C^m$ (c.f.r. Section~\ref{sec:orbitrep}). In particular we can suppose that  
\[f_h(v_j)= \begin{cases} 0, \mbox{ if } j\in \{0, \dots k\}; \\ v_j, \mbox{ otherwise.}\end{cases}\]
Set $B_1=\{v_1, \dots, v_{d_1}\}$ and define recursively 
\[ B_{i+1} = \begin{cases}   B_{i} \cup \{v_{d_i +1}, \dots, v_{d_{i +1}} \, v_m\} \setminus \{v_1\}, &\mbox { if } i = h+1; \\ B_{i} \cup\{v_{d_i +1}, \dots, v_{d_{i +1}}\}, &\mbox{ otherwise. } 
\end{cases}\]
Consider the $n$-tuple of subspaces $F=(F_1, \dots, F_n)$ such that $F_i=\Span B_i$.
By definition, the indecoposable representation $U_{1,h}$ appears as an indecoposable direct summand in the subrepresentation $L_{F}$ of $M_\mathbf{f}$. Moreover, since $v_1 \notin F_{h+1}$, then there exists $n'$ such that $U_{h+1,n'}$ is an indecomposable direct summand of $M_\mathbf{f} / L_F$. Since $\mrm{Ext}^1 (U_{1,h} , U_{h+1,n'})  \neq 0$, we have that $\mrm{Ext}^1 (L_F , M_{\mathbf{f}}/L_F)  \neq 0$ and $\mathcal{F}l_\mathbf{d}^\mathbf{f} (\C^m)$ is not smooth at $F$.
\endproof
\begin{corollary}A linear degeneration $\mathcal{F}l_\mathbf{d}^\mathbf{f} (\C^m)$ is smooth if and only if it is a product of partial flag varieties.
\end{corollary}
\begin{rmk} In the proof of the Theorem~\ref{thm:smooth} we proved a stronger statement: whenever $\rk f_i < m$ for some $i$, there exists a point $L \in Gr_{\mathbf{d}}(M_\mathbf{f})$ such that $\mrm{Ext}^1 (L , M_{\mathbf{f}}/L)  \neq 0$. Observe that in general this does not implies that the point $L$ is singular in $Gr_{\mathbf{d}}(M_\mathbf{f})$. 
\end{rmk}

\section{Stratification of Space of Representations}\label{sec:strata}
As a result of the classification of smooth linear degenerations, we describe a new stratification of the total space of representations $R:= \left(\mrm{End}(\C^m)\right)^{n-1}$.
Let $I$ be a subset of $\{1, \dots, n-1\}$. For each $I$ we can associate an $n-1$-tuple of endomorphisms $\mathbf{f}_I = (f_1, \dots, f_{n-1})$ such that 
\[ f_i = \begin{cases} 0, &\mbox{ if }  i \in I;\\ \mrm{id}, &\mbox{ otherwise.} \end{cases}\]
Once a dimension vector $\bf{d}$ is fixed, the linear degeneration $\mathcal{F}l_\mathbf{d}^{\mathbf{f}_I} (\C^m)$ is smooth, as shown by Theorem~\ref{thm:smooth}. We denote by $O_I$ the orbit of $\mathbf{f}_I$ in $R$. Theorem~\ref{thm:smooth} further implies that every smooth linear degeneration is the fiber of an $n-1$-tuple of endomorphisms $\bf f$ contained in some orbit of the form $O_I$.
For each $I$, we define $S_I$ as the subset of $R$ given by 
\[S_I=\{ (f_1, \dots, f_{n-1})\in R \, | \, f_i = 0 \Longleftrightarrow i \in I \}.\]
Observe that each $S_I$ is a locally closed subset of $R$ and that $\cup_I S_I = R$. Furthermore, since $S_I \subset \overline{S_J}$ if and only if $J \subset I$, the family $\{S_I\}_{I \subset \{1, \dots n-1\}}$ provides a stratification for $R$. The poset of such strata, ordered by closure inclusion, is isomorphic to the reverse Boolean poset over $n-1$ elements. Additionally, each stratum is $G$-invariant, and $O_I$ is the unique dense open orbit in $\overline{S_I}$.

 The flat locus $U^I_{flat}$ in $\overline{S_I}$ is the set of $\mathbf{f} \in \overline{S_I}$ such that $\mathcal{F}l_\mathbf{d}^\mathbf{f} (\C^m)$ has minimal dimension. The flat irreducible locus $U^I_{flat, irr}$ is the subset of $U^I_{flat}$ such that $\mathcal{F}l_\mathbf{d}^\mathbf{f} (\C^m)$ is also an irreducible variety.
If $I=\{i_1 < \dots < i_k\}$, the set $\overline{S_I}$ can be identified with the subspace of $R$ given by the product:
\[R_I:=End(\C^m)^{i_1-1} \times End(\C^m)^{i_2-i_1-1}\times \dots \times End(\C^m)^{n-i_k-1}.\]
 Set $Y_I=\pi^{-1}(\overline{S_I})$, the universal linear degeneration map $\pi$ restricts to a morphism of algebraic varieties 
\[\pi_I: Y_I \rightarrow R_I.\]
By definition, this morphism is $G$-equivariant, and since $Y_I$ is a closed subset of $Y$, it is also a proper map. Set $\mathbf{d}_{j+1}=(d_{i_j+1} \dots d_{i_{j+1}})$, where $j \in \{0, \dots , k\}$ and $i_0=0$. For each $j$, consider the universal linear degeneration associated with the dimension vector $\mathbf{d}_j$:
\[ \pi_{j}: Y_{j} \rightarrow R_j:=\mrm{End}(\C^m)^{i_j-i_{j-1}-1}.\]
As a consequence of Remark~\ref{rmk:product}, we have: 
\[Y_I= Y_1 \times \dots \times Y_k, \qquad \pi_I = \pi_1 \times \dots \times \pi_k.\]
In particular, $Y_I$ is an irreducible, smooth and Cohen-Macaulay variety. The flat and flat irreducible loci described in \cite{CIFFFR20} and \cite{CIFFFR17} appear as special cases of flat and flat irreducible loci in $\overline{S_\emptyset}$. 

Furthermore, denote by $U_{flat}^j$ and $U_{flat, irr}^j$ the flat and the flat irreducible loci in $R_j$, respectively.  We have that
\[U^I_{flat} = U_{flat}^1 \times \dots \times U_{flat}^k,  \qquad U^I_{flat, irr} = U_{flat, irr}^1 \times \dots \times U_{flat, irr}^k, \]and $\pi_I$ is flat when restricted to~$U^I_{flat}$.

According to Theorem~\ref{thm:flatlocus}, there exist two orbits $O_{r^2}^j, O_{r^1}^j \subset R_j$ such that $U_{flat}^j$ is the locus of orbits degenerating to $O_{r^2}^j$, and $U_{flat, irr}^j$ is the locus of orbits degenerating to  $O_{r^1}^j$.

We recall that $U_{flat}$ and $U_{irr, \, flat}$ are characterized in \cite{CIFFFR20} and \cite{CIFFFR17} also using rank sequences. An analogous description is also possible in this more general context. 
Consider the following rank series: 


\[r^{1,I}_{i}=m-d_{i+1}+d_i, \qquad r^{1,I}_{i,j}=\begin{cases} m - d_{j+1}+ d_i, &\mbox{ for } i_h+1 \leq  i < j  < i_{h+1};  \\ 
0, &\mbox{ otherwise, }\end{cases}\]

\[r^{2,I}_{i}=m-d_{i+1}+d_i -1, \qquad r^{2,I}_{i,j}=\begin{cases} m - d_{j+1} + d_i -1, &\mbox{ for } i_h+1 \leq  i < j  < i_{h+1};  \\ 
0, &\mbox{ otherwise. }\end{cases}\]
Observe that the rank sequences $r^{1,I}$ and $r^{2,I}$ correspond to the orbits 
\[O_{r^{1,I}}= O_{r^{1}}^1 \times \dots \times O_{r^{1}}^k, \qquad O_{r^{2,I}}= O_{r^{2}}^1 \times \dots \times O_{r^{2}}^k, \]
respectively. As a consequence, we obtain the following:
\begin{theorem}
\begin{enumerate}
\item The flat locus in $\overline{S_I}$ is the union of orbits $O_\mathbf{r}$ degenerating to $O_{\mathbf{r}^{2,I}}$,
\item The flat irreducible locus in $\overline{S_I}$ is the union of orbits $O_\mathbf{r}$ degenerating to $O_{\mathbf{r}^{1,I}}$.
\end{enumerate}
\end{theorem}
Note that every irreducible linear degeneration appears as a fiber of certain ${\bf f} \in U^I_{flat, irr}$ for a suitable~$I$. In particular, the subset 
\[ U_{irr}:= \bigcup_{I \subset \{1, \dots, n-1\}} U^I_{flat, irr} \subset R \]
coincides with the locus of tuples $\bf f$ such that $ \mathcal{F}l_\mathbf{d}^\mathbf{f} (\C^m)$ is irreducible.

\section{Regularity in the Irreducible Locus}
An algebraic variety $X$ is said to be regular in codimension $k$ if there exists a codimension $k+1$ closed subvariety $X' \subset X$ such that the singular locus of $X$ is contained in $X'$. In this section, we discuss the regularity of linear degenerations $\mathcal{F}l_\mathbf{d}^\mathbf{f} (\C^m)$ when $\mathbf{f}$ lies in the irreducible locus $U_{irr}$. 
\subsection{Preliminaries on Morphisms of Schemes}
We firstly recall some scheme theoretic results, which we are going to use to describe the singular locus in flat irreducible linear degenerations.
\begin{lemma}\cite[\href{https://stacks.math.columbia.edu/tag/05F6}{Section 05F6}]{stacks-project}\label{lem:upsemicontinuous} Let $f: X \rightarrow Y$ be a proper morphism of schemes. The function $\eta : Y \rightarrow \mathbb{N}$ that associates to $y \in Y$ the dimension of the fiber $X_y$ is upper semi-continuous.
\end{lemma}
In particular, for any positive integer $d$, the set $Z_d = \{y \in Y \; | \; \dim X_{y} > d\}$ is a closed subsets in~$Y$. 
\begin{proposition}\label{prop:dimorbite} Let $G$ be a complex algebraic group and $X, Y$ be two complex $G$-varieties. Let $f:X \rightarrow Y$ be a proper $G$-equivariant morphism. Let $y, y' \in Y$ and suppose that the $G$-orbit $O_y$ is contained in the closure of $O_{y'}$. Then $\dim X_y \geq \dim X_{y'}$. 
\end{proposition}
\proof 
 Suppose $d=\dim X_y < \dim X_{y'}$ and consider set $Z_d$ as above. Since $f$ is $G$-equivariant, $Z_d$ is a $G$-invariant subset of $Y$. By definition, we have that ${y'} \in Z_d$ and consequently  $O_{y'} \subset Z_d$ by $G$ invariance.  Since $Z_d$ is closed, the orbit closure $\overline{O_{y'}}$ is contained in $Z_d$. By hypothesis $O_y \subset \overline{O_{y'}}$ and then $O_y \subset Z_d$, in contradiction with our assumption. 
\endproof 
\begin{theorem}\cite[Theorem 12.2.4]{EGA}\label{th:grazieGrot} Let $f : X \rightarrow Y$ be a flat proper morphism of varieties. Then the locus of $y \in Y$ such that $f^{-1}(y)$ is reduced (resp. irreducible, resp. normal) is open in~$Y$. 
\end{theorem}
\begin{proposition}\label{prop:reducedfibers} Let $G$ be a complex algebraic group and $X, Y$ two complex $G$-varieties. Let $f:X \rightarrow Y$ be a flat proper $G$-equivariant morphism. Suppose that there exists $\overline{y} \in Y$ such that: 
\begin{itemize}
\item the fiber $f^{-1}(\overline{y})$ is a reduced (resp. irreducible, resp. normal) subvariety of $X$,
\item the orbit $O_{\overline y}$ is is contained in the closure of every $G$-orbit in $Y$,
\end{itemize}
then $f^{-1}(y)$ is reduced (resp. irreducible, resp. normal) for every $y$ in $Y$.
\end{proposition}
\proof 
Let $V \subset Y$ be the locus of points $y \in Y$ such that $f^{-1}(y)$ is reduced (resp. irreducible, resp. normal).
Since $f$ is $G$-equivariant, $V$ is a $G$-invariant subset and $O_{\overline{y}} \subset V$. 
Because of Theorem~\ref{th:grazieGrot} it is an open subset of $Y$, then its complementary $V^c$ is a $G$-invariant closed subset of $Y$. 
Suppose that $V^c$ is not empty, since it is closed, it contains the closure of orbits $O_{y'}$ for any $y' \in V^c$. But this implies $O_{\overline y} \subset V^c$, which is absurd.
\endproof 

We recall now some results about the singular locus of morphisms of schemes.

\begin{defn}
Let $Y$ be a locally Notherian scheme. A morphism of schemes $f:X \rightarrow Y$ of locally finite type is smooth at $x\in X$ if: \begin{enumerate} \item  the local map $O_{Y,f(x)} \rightarrow O_{X,x}$ is flat, \item the $k(f(x))$ scheme $X_{f(x)} \rightarrow k(f(x))$ is smooth at $x$. \end{enumerate} \end{defn} 
In particular if $f$ is a flat morphism between varieties over an algebrically closed field, smoothness of $f$ at $x$ is equivalent to prove that $X_{f(x)}$ is a smooth variety at $x$ (c.f.r. \cite[Chapter 4 - Definition 3.35]{Liu}). Since over algebrically closed fields notions of regularity and of smoothness coincide (see \cite[Chapter 4 - Corollary 2.17]{Liu}), a flat morphism $f$ is smooth if and only if for every $y \in Y$ the fiber $X_y$ is smooth at its closed points.
\begin{defn}
Let $f : X \rightarrow Y$ be a locally of finite type morphism of schemes. The singular locus of $f$ is the set 
\[\mrm{Sing}(f)=\{ x \in X \;|\;  f \mbox{ is not smooth at } x \}.\] 
\end{defn}
It is a classical fact that $\mrm{Sing}(f)$ is a closed subset of $X$ (see \cite[\href{https://stacks.math.columbia.edu/tag/01V5}{Definition 01V5}]{stacks-project}).
\subsection{The Singular Locus in Well Behaved Linear Degenerations}
We now recall some key results concerning the regularity of linear degenerations within the orbit $O_{\mathbf{r}^1}$. These degenerations can be identified with "well behaved" quiver Grassmannians, studied in \cite{CIFR12} and \cite{CIFR13}. 
\begin{defn}\label{def:Feigin} Let $\mathbf{d}=(d_1< \dots< d_n)$ be a dimension vector, where $d_n < m$. A linear degeneration $\mathcal{F}l_\mathbf{d}^{\mathbf{f}} (\C^m)$ is well-behaved if $\rk f_i = m-(d_{i+1}-d_i)$ for every $i$ and the kernels of the maps $f_i$ are linearly independent. 
\end{defn}
For a well-behaved linear degeneration, the associated representation $M_\mathbf{f}$ can be decomposed as $M_{\bf f}=P
\oplus I$, where $P$ is a projective subrepresentation and $I$ is an injective subrepresentation. Specifically, we have: \[M_\mathbf{f}= \left( \bigoplus_{i=1}^{n-1} U_{1,i}^{\oplus d_{i+1}-d_i} \right) \oplus U_{1,n}^{m-d_n+d_1} \oplus \left(\bigoplus_{i=1}^{n-1} U_{i+1,n}^{ \oplus d_{i+1}-d_i} \right) \]
and \[I= U_{1,n}^{d_1} \oplus \bigoplus_{i=1}^{n-1} U_{1,i}^{\oplus d_{i+1}-d_i}, \qquad \qquad P=U_{1,n}^{m-d_n } \oplus \bigoplus_{i=1}^{n-1} U_{i+1,n}^{ \oplus d_{i+1}-d_i}.  \]
A well behaved linear degeneration $\mathcal{F}l_\mathbf{d}^{\mathbf{f}} (\C^m)$ can then be identified with the quiver Grassmannian $\mrm{Gr}_{\dim P}(P \oplus I)$. 
These quiver Grassmannians are known as "well behaved" (c.f.r. \cite[Section 3 ]{CIFR12}) and their properties have been extensively investigated in~\cite{CIFR12}, in the more general context of representations of Dynkin quivers. 
Note that the representation $M_\mathbf{f}$ is catenoid, which implies that a well behaved quiver Grassmannian is a Schubert quiver Grassmannian. We now summarize several key properties of these quiver Grassmannians:
\begin{proposition}\cite[Proposition 3.1]{CIFR12}\label{prop:P+I}
Let $M$ be a representation such $M= P \oplus I$, where $P$ is projective and $I$ is injective. 
\begin{enumerate} 
\item  The dimension of $\mrm{Gr}_{\bf dim P}(M)$ equals $\langle \mathbf{dim P} , {\bf{dim M}}- \bf{dim P} \rangle $,
\item $\mrm{Gr}_{\bf{dim P}}(M)$ is reduced, irreducible and rational,
\item  $\mrm{Gr}_{\bf{dim P}}(M)$ is a locally complete intersection scheme.
\end{enumerate}
\end{proposition}
\begin{proposition}\cite[Theorem 4.5]{CIFR12}\label{prop:normal}
A quiver Grassmannian of the form $\mrm{Gr}_{\bf dim P}(P \oplus I)$ is regular in codimension 1 and consequently it is a normal variety. 
\end{proposition}
The statement of \cite[Theorem 4.5]{CIFR12} was further generalized in~\cite{CIFR13}:
\begin{theorem}\cite[Theorem A.1]{CIFR13}\label{thm:reg2} A quiver Grassmannian of the form $\mrm{Gr}_{\bf dim P}(P \oplus I)$ is regular in codimension 2.\end{theorem}
More precisely, the authors prove that if $N$ is a singular point in $\mrm{Gr}_{\bf dim P}(P \oplus I)$, then the closure of its isostratum is a subvariety of codimension at least 3. 
\subsection{Normality of Irreducible Linear Degenerations}
Since well behaved degenerations precisely correspond to $n-1$-tuples in $O_{\bf r}^1$, as a consequence of Proposition~\ref{prop:reducedfibers} we have the following: 
\begin{theorem}\label{thm:reducedegenerations}
For any ${\bf f} \in U_{flat, irr}$ the linear degeneration $\mathcal{F}l_\mathbf{d}^{\mathbf{f}} (\C^m)$ is normal.
\end{theorem}
The previous Theorem is essentially a generalization of~\cite[Theorem B]{CIFFFR17} and its proof is provided here for completeness.
\proof 
Since the property of a morphism of being proper is stable under base change, the restriction \[ \pi : \pi^{-1}(U_{flat, irr}) \rightarrow U_{flat, irr}\] is a flat proper $G$-equivariant morphism between complex $G$-varieties. Consider $ \overline{\bf f} \in O_{\bf r^1}$. The linear degeneration $\mathcal{F}l_\mathbf{d}^{\overline{\mathbf{f}}} (\C^m)$ is a well behaved degeneration and it is normal by Proposition~\ref{prop:normal}. Because of~\cite[Theorem A]{CIFFFR20}, $ O_{\bf r^1}$ is contained in the closure of every orbit in $U_{flat, irr}$. We are under the hypotheses of Proposition~\ref{prop:reducedfibers} and the theorem follows.
\endproof
This theorem has remarkable consequences regarding the normality of irreducible linear degenerations.
\begin{theorem}Every irreducible linear degeneration is a normal variety. 
\end{theorem}
\proof 
According to Corollary~\ref{cor:euler}, every irreducible linear degeneration is a product of flat irreducible linear degenerations. Since the product of normal varieties is normal, the statement follows from Theorem~\ref{thm:reducedegenerations}.  \endproof
\subsection{The Singular Locus of Certain Quiver Grassmannians}
Let $X$ be an algebraic variety, we denote its singular locus by $\mrm{Sing}(X)$. We recall that, once a flat irreducible degeneration $\mathcal{F}l_\mathbf{d}^\mathbf{f} (\C^m)$ is identified with the quiver Grassmannian $\mrm{Gr}_{\mathbf{d}}(M_{\mathbf{f}})$ as in Section~\ref{sec:smooth}, a point $N$ is singular  
if and only if $\mathrm{Ext}^1\left(N, M_\mathbf{f}/N\right)\neq 0$. Consequently, if ${\bf{f}} \in U_{flat, irr}$, we have that
\[\mrm{Sing} \left(\mrm{Gr}_{\mathbf{d}}(M_{\mathbf{f}}) \right) = \{ N \in \mrm{Gr}_{\mathbf{d}}(M_{\mathbf{f}}) |  \mathrm{Ext}^1\left(N, M_\mathbf{f}/N\right)\neq 0 \}\]
For certain quiver Grassmannians, the singular locus can be explicitly described. Fix a basis $\{v_1, \dots, v_m\}$ of $\C^m$ and consider $h \in \{1, \dots, n-1 \}$. Let $\mathbf{f}^h$ be the $n-1$-tuple of endomorphisms such that 
\[f^h_i=
\begin{cases}
\pi_1, &\mbox{ if } i = h; \\
\mrm{id}, &\mbox{ otherwise.}
\end{cases}\]
where $\pi_1$ is the projection along the vector $v_1$. For brevity, we  denote the representation $M_{\mathbf{f}^h}$ by $M^h$. Observe that $M^h$ decomposes as 
\[M^h= U_{1,n}^{m-1} \oplus U_{1,h} \oplus U_{h+1,n}.\] 
and, in particular, it is a catenoid representation.
As a consequence, given a dimension vector $\mathbf{d}=(d_1 < \dots < d_n)$, the quiver Grassmannian $\mrm{Gr}_{\mathbf{d}}(M^h)$ is a Schubert quiver Grassmannian. Moreover, because of classification provided in~\cite{CIFFFR20}, it can be realized as a flat irreducible linear degeneration of the partial flag variety~$\mathcal{F}l_{\mathbf{d}}(\C^m)$.
\begin{rmk}\label{rmk:upperdeg} Let $X=\mathcal{F}l_\mathbf{d}^\mathbf{f} (\C^m)$ be a linear degeneration such that $\bf f$ is not in the orbit of $(\mrm{id}, \dots, \mrm{id})$. By examining rank sequences, we can verify that if $\rk f_h < m$, then $X$ is a degeneration of  $\mrm{Gr}_{\mathbf{d}}(M^h)$. 
\end{rmk}
We will now explicitly describe the singular locus of these specific degenerations. For this purpose, we first present some preliminary lemmas.  The first lemma holds generally for representations of equioriented quivers of type $A_n$ and follows from basic notions of quiver representation theory.
\begin{lemma}\label{lem:surj} Consider a representation $M$ of the equioriented quiver of type $A_n$. The representation $M$ does not have an indecomposable direct summand of the form $U_{k+1,b}$ for $b \in \{k+1, \dots, n\}$ if and only if the map $M_{k,k+1}$ is surjective.
\end{lemma}
The second one holds true only for representations of type $M^h$. Denote by $M_k$ the vector space associated to the vertex $k$ in $M^h$. For each $M_k$ in $M^h$ fix now a basis $\{v^k_1, \dots, v^k_m\}$. We denote by $\rho$ the projection operator $$\rho:M_{h+1} \rightarrow  \Span \{ v^{h+1}_1 \} $$ defined by 
\[\rho(v^{h+1}_j)=
\begin{cases}
v^{h+1}_j, &\mbox{ if } j = 1; \\
0, &\mbox{ otherwise};
\end{cases}\]
on the basis $\{v^{h+1}_1, \dots, v^{h+1}_m\}$ and extended by linearity to the whole $M_{h+1}$. 
\begin{lemma}\label{lem:projection} Let $N$ be a subrepresentation of $M^h$ with dimension vector $\mathbf{d}=(d_1< \dots <d_n )$. Assuming that $v^{h}_1 \in N_h$, then there exists $b \in \{h+1, \dots, n\}$ such that the representation $U_{h+1,b}$ is an indecomponsable direct summand of $M/N$ if and only if the restriction of $\rho$ to $N_{h+1}$ is the null map.
\end{lemma}
\proof Suppose firstly that the restriction of $\rho$ to $N_{h+1}$ is the null map, i.e. $N_{h+1} \subset \Span \{ v^{h+1}_2, \dots v^{h+1}_m\} $. In virtue of Lemma~\ref{lem:surj}, it is enough to prove that the map 
\[\widetilde{f_h}: M_h/N_h \longrightarrow M_{h+1}/N_{h+1}\]
induced by $f_h$ on quotient is not surjective. 
Without losing of generality, we can suppose $N_h= \mrm{Span} \{ v_1^h,u_1, \dots, u_{d_h-1}\}$ with $u_i  \in \Span \{ v^{h}_2, \dots v^{h}_m\}$ for every $i$. 
Set $w_i=f_h(u_i)$.
By definition of $f_h$ we can also suppose that
 $$N_{h+1}= \Span \{ w_1, \dots , w_{d_h-1}, \overline{w}_{d_h} \dots \overline{w}_{d_{h+1}}\}$$ 
 for certain vectors $\overline{w}_{d_h} \dots \overline{w}_{d_{h+1}}$.
Since $\rho$ restricted to $N_{h+1}$ is the null map, we have that $$\{ w_1, \dots , w_{d_h-1}, \overline{w}_{d_h} \dots \overline{w}_{d_{h+1}}\} \subset \Span \{ v^{h}_2, \dots v^{h}_m\}.$$
We can complete the basis of $N_{h+1}$ to a basis of $M_{h+1}$ choosing $m-d_{h+1} -1$ vectors $ \widetilde{w}_1, \dots, \widetilde{w}_{m-d_{h+1}-1} \in  \mrm{Span}\{ v^{h+1}_2, \dots v^{h+1}_m\}$ such that $$M_{h+1}= N_{h+1} \oplus \Span \{ v^{h+1}_1\} \oplus \Span \{ \widetilde{w}_1, \dots, \widetilde{w}_{m-d_{h+1}-1}\}.$$ 
 By definition of $f_h$, there exist $\overline{u_1}, \dots, \overline{u}_{m-d_{h+1}-1} \in \Span \{ v^{h}_2, \dots v^{h}_m\} $ such that $f_h(\overline{u}_i)=\widetilde{w}_i$ and $$\{ v_1^h,u_1, \dots, u_{d_h-1}, \overline{u}_1, \dots, \overline{u}_{m-d_{h+1}-1} \} $$ is a set of linearly independent vectors in $M_h$. 
 We can finally complete this set to a basis 
 \[B_h = \{ v_1^h,u_1, \dots, u_{d_h-1}, \overline{u_1}, \dots, \overline{u}_{m-d_{h+1}-1} \} \cup \{ \widetilde{u}_1, \dots, \widetilde{u}_{d_{h+1}-d_h}\}\] for $M_h$, where again $ \{ \widetilde{u}_1, \dots, \widetilde{u}_{d_{h+1}-d_h}\} \subset \Span \{ v^{h}_2, \dots v^{h}_m\}$.
Passing to the quotient, we have that $[\overline{u_1}], \dots, [\overline{u}_{m-d_{h+1}-1}] , [\widetilde{u}_1], \dots, [\widetilde{u}_{d_{h+1}-d_h}]\}$ and $\{ [v^{h+1}_1], [\widetilde{w}_1], \dots, [ \widetilde{w}_{m-d_{h+1}-1}]\}$ are basis of $M_h/N_h$ and $M_{h+1}/N_{h+1}$ respectively, where by $[v]$ we denote the class of $v$ in the quotient. By construction we have that $[v^{k+1}_1] \notin \mrm{Im} \widetilde{f_h}$, so $f_h$ is not surjective and there exists $b \in \{h+1, \dots, n\}$ such that $U_{h+1, b}$ is a subrepresentation of $M/N$. \\
Conversely, suppose that for some $b \in \{k+1, \dots, n\}$ the representation $U_{k+1, b}$ is an indecomposable summand of $M/N$ and that $N_{h+1}$ is not contained in $\mrm{Ker}\, \rho$. In particular, there exists a vector $u \in N_{h+1}$ such that $\rho(u)=v^{h+1}_1$. Because of definition of $\rho$, we can write $u=v^{h+1}_1 + w$, with $w \in \Span \{ v^{h+1}_2, \dots v^{h+1}_m\}$. Choose now $d_h-1$ linearly independent vectors $v_1, \dots, v_{d_h-1}$ in $\mrm{Span}\{v^h_2, \dots, v^h_m\}$ such that $N_{h}=\Span \{ v^{k}_1, v_1, \dots, v_{d_h-1} \} $. 
Without loss of generality, we can also suppose $$N_{h+1}= \Span \{ u, w_1, \dots, w_{d_{h}-1}, \overline{w}_{d_h+1}, \dots, \overline{w}_{d_{h+1}} \},$$ where $w_i=f_h(v_i)$ for every $i\leq d_h-1$.
Observe that, by definition of $f_h$, we have that $w_i \in \Span \{ v^{h+1}_2, \dots v^{h+1}_m\}$. Furthermore, by Gauss elimination, we can suppose that $\overline{w}_i \in \Span \{ v^{h+1}_2, \dots v^{h+1}_m\}$. We choose now $m-d_{h+1}$ linearly independent vectors $\{u_1, \dots u_{m-d_{h+1}} \}$ such that $$B_{h+1}= \{ u, w_1, \dots, w_{d_{h}-1}, \overline{w}_{d_h+1}, \dots, \overline{w}_{d_{h+1}} \} \cup \{u_1, \dots u_{m-d_{h+1}} \}$$ is a basis for $M_{h+1}$.
Since $v^{h+1}_1 + w \in N_{h+1}$, again by Gauss elimination, we can choose $\{u_1, \dots u_{m-d_{h+1}} \}$ in $\Span \{ v^{h+1}_2, \dots v^{h+1}_m\}$. 
Since $f_h$ is the identity when restricted to $\mrm{Span}\{v^h_2, \dots, v^h_m\}$, there exists two sets of linearly independent vectors $\{\widetilde{w}_{d_h+1}, \dots, \widetilde{w}_{d_{h+1}} \}$ and $\{\overline{u}_1, \dots \overline{u}_{m-d_{h+1}} \}$ such that: 
\begin{itemize}
\item $\overline{u_i} \notin N_h$ and $f_h(\overline{u}_i)=u_i$ for every $i \in \{1, \dots m-d_{h+1} \}$,
\item $\widetilde{w}_j \notin N_h$ and $f_h(\widetilde{w_j})=\overline{w}_j$ for every $j \in \{d_h+1, \dots d_{h+1} \}$.
\end{itemize}
Consequently we can complete $\{ v^{h}_1, v_1, \dots, v_{d_h-1}\}$ to a basis $B$ of $M_h$ choosing 
\[B_h=\{ v^{h}_1, v_1, \dots, v_{d_h-1}\} \cup \{\widetilde{w}_{d_h+1}, \dots, \widetilde{w}_{d_{h+1}} \} \cup \{\overline{u}_1, \dots \overline{u}_{m-d_{h+1}} \}. \]
and passing to the quotient, it results immediately that the map $\widetilde{f}_h$ induced by $f_h$ is surjective. Since we supposed $U_{h+1,b}$ to be an indecomposable summand of $M/N$, this is in contrast with Lemma~\ref{lem:surj} and consequently $N_{h+1}$ need to be contained in $\mrm{Ker} \,\rho$.
\endproof
Consider now the subrepresentation $M' \subset M^h$ defined by 
\[ M'_i= \begin{cases} \Span \{ v^i_2, \dots, v^i_m \}, \mbox{ if }  i=h,h+1; \\ M_i, \mbox{ otherwise.}\end{cases} \]
We have
$$M' = U_{1,n}^{m-1} \oplus U_{1,h-1} \oplus U_{h+2,n}.$$  Moreover, consider the dimension vector $\mathbf{d'}$ defined by 
\[d'_i=\begin{cases}d_i-1 &\mbox{ if } i=h \\ d_i &\mbox{ otherwise.}\end{cases} \]
\begin{theorem}
The singular locus $\mrm{Sing}(\mrm{Gr}_{\mathbf{d}}(M^h)) $ is a closed subvariety, isomorphic to the quiver Grassmannian $\mrm{Gr}_{\mathbf{d'}}(M')$.
\end{theorem}
\proof
Let $N$ be a point in $\mrm{Gr}_{\mathbf{d}}(M^h)$. Since $\mrm{Gr}_{\mathbf{d}}(M^h)$ can be identified with a flat irreducible linear degeneration, $N$ is singular if and only if $\mathrm{Ext}^1\left(N, M^h/N\right)\neq 0$.
Since $$M^h= U_{1,n}^{m-1} \oplus U_{1,h} \oplus U_{h+1,n},$$ the condition $\mathrm{Ext}^1\left(N, M^h/N\right)\neq 0$ is satisfied if and only if $U_{a,h}$ and $U_{h+1,b}$ are indecomposable components of $N$ and $M^h/N$ respectively, for some positive integers $a \in \{1, \dots h \}$ and $ b \in \{h+1,\dots, n\}$.
There exists an $a$ such that indecomposable representation $U_{a,h}$ is a component of $N$ if and only if $v^h_1 \in N_h$. Moreover, because of Lemma~\ref{lem:projection}, there exists a $b$ such that $U_{h+1,b}$ is an indecomposable component of $M^h/N$ if and only if $\rho: N_{h+1} \rightarrow \mrm{Span} \{ v^{h+1}_1 \} $ is the null map. 
Therefore, we have 
\[\mrm{Sing}(\mrm{Gr}_{\mathbf{d}}(M^h)) = \left\{N \subset M \; | \mathbf{dim N} = \mathbf{d}, \; \dim \left( N_h \cap \Span \{ v^h_1\}\right) = 1, \, N_{h+1}  \subset Ker(\rho) \right\}\]
If $N' \in M'$, define \[\Sigma(N')_i=\begin{cases}N'_i \oplus \mrm{Span}(v_1^k), &\mbox{ if } i = h; \\ N'_i, &\mbox{ otherwise.} \end{cases}\] where $N'_h$ and $N'_{h+1}$ are identified with subspaces of $M_{h}$ and $M_{h+1}$ respectively, via the immersion of $\Span \{ v^h_2, \dots, v^h_m \}$ into $\Span \{ v^h_1, \dots, v^h_m \}$. Consider now the map
\[\begin{matrix}
\Sigma:& \mrm{Gr}_{\mathbf{d'}}(M') &{\longrightarrow} &\mrm{Gr}_{\mathbf{d}}(M_{\mathbf{f}}^h),\\
&N'=(N'_i) &\longrightarrow & (\Sigma(N')_i).
\end{matrix} \]
We have that $\Sigma$ defines a morphism of varieties and that its image is contained in $\mrm{Sing}(\mrm{Gr}_{\mathbf{d}}(M^h))$.
The inverse morphism is provided by considering the map 
\[\begin{matrix}
\Sigma':& \mrm{Sing} \left( \mrm{Gr}_{\mathbf{d}}(M_{\mathbf{f}}^h)\right) &{\longrightarrow} &\mrm{Gr}_{\mathbf{d}'}(M'),\\
&(N_1, \dots, N_n) &\longrightarrow &(N_1, \dots, N_{h-1}, \pi_1(N_h), N_{h+1}, \dots, N_n),
\end{matrix} \]
and the desired isomorphism is proved.
\endproof
\begin{corollary}\label{cor:codimsing}
The singular locus $\mrm{Sing}(\mrm{Gr}_{\mathbf{d}}(M_{\mathbf{f}}^h))$ is an irreducible smooth closed subvariety of $\mrm{Gr}_{\mathbf{d}}(M_{\mathbf{f}})$ of codimension $2(d_{h+1}-d_h)+1$. In particular, if $d_{h+1}-d_h=1$, then $\mrm{Sing}(\mrm{Gr}_{\mathbf{d}}(M_{\mathbf{f}}^h))$ has codimension~3.
\end{corollary}
\proof Since the representation $M'$ is rigid, the quiver Grassmannian $\mrm{Gr}_{\mathbf{d}'}(M')$ is irreducible and smooth by \cite[Proposition 2.3]{CI20}. Let us denote by $\mathbf{\epsilon}_j$ the element of $\Z^n$ with all coordinates equal to zero except for the $j$-th one, that is equal to 1. Set $\mathbf{m}:=\mathbf{dimM}=(m, \dots, m)$ and $\mathbf{dim M'}=\mathbf{m}'=\mathbf{m}- \mathbf{\epsilon}_h - \mathbf{\epsilon}_{h+1}$. The dimension of the quiver Grassmannian $\mrm{Gr}_{\mathbf{d'}}(M')$ equals $\langle \mathbf{d}', \mathbf{m}' - \mathbf{d}'\rangle $. Consequently we have 
\begin{align*} \langle \mathbf{d}', \mathbf{m}' - \mathbf{d}'\rangle &= \langle \mathbf{d}', \mathbf{m} - \mathbf{d}' - \mathbf{\epsilon}_h - \mathbf{\epsilon}_{h+1} \rangle \\
&= \langle \mathbf{d}-\mathbf{\epsilon}_{h} , \mathbf{m} - \mathbf{d}  - \mathbf{\epsilon}_{h+1} \rangle \\
&= \langle \mathbf{d} , \mathbf{m} - \mathbf{d} \rangle  - \langle \mathbf{d}, \mathbf{\epsilon}_{h+1} \rangle  -  \langle \mathbf{\epsilon}_{h} , \mathbf{m} - \mathbf{d}  - \mathbf{\epsilon}_{h+1} \rangle \\
&= \langle \mathbf{d} , \mathbf{m} - \mathbf{d} \rangle  - \langle \mathbf{d}, \mathbf{\epsilon}_{h+1} \rangle  - \left[ \langle \mathbf{\epsilon}_{h} , \mathbf{m} \rangle - \langle \mathbf{\epsilon}_{h} ,\mathbf{d} \rangle  - \langle \mathbf{\epsilon}_{h} ,\mathbf{\epsilon}_{h+1} \rangle \right] \\
&= \langle \mathbf{d} , \mathbf{m} - \mathbf{d} \rangle  - 2(d_{h+1}-d_h) -1.
\end{align*}
The statement now follows, recalling that $\mrm{Gr}_{\mathbf{d}}(M_{\mathbf{f}}^h)$ is a flat irreducible linear degeneration of $\mathcal{F}l_{\bf d}(\C^m)$ and its dimension equals $\langle \mathbf{d} , \mathbf{m} - \mathbf{d} \rangle$. 
\endproof
\subsection{The Singular Locus of Flat Irreducible Degenerations}\label{subsec:regularity}
We can now prove the main theorems of this section. Let $Y$ be the universal linear degeneration, recall that the map 
\[\pi: Y \rightarrow R:= \left( End (\C^m) \right)^{n-1}\]
is a proper $G$-equivariant morphism of complex algebraic varieties. 
\begin{proposition}\label{prop:closuredim}
Let $X, X'$ be two flat linear degenerations. If $X$ degenerates to $X'$ then
\[\dim \mrm{Sing}(X') \geq \dim \mrm{Sing}(X).\]
\end{proposition}
 \proof Let us denote by $Z$ the singular locus of $\pi$ and let $i$ be the closed immersion of $Z$ into $Y$. Since closed immersions are proper morphisms and composition of proper morphisms is proper, the map 
 $\pi \circ i: Z \rightarrow Y \rightarrow R$
 is a proper morphism of schemes. Moverover, since $\pi$ is $G$-equivariant, then $\pi \circ i$ is $G$-equivariant too.
 Consider now $\bf f$ and $\bf f ' $ such that $O_\mathbf{f} \subset \overline{O_{\mathbf{f}'}}$. Since we are in the hypotheses of Proposition~\ref{prop:dimorbite}, then \begin{equation}\label{eq:dimsing} \dim (\pi\circ i)^{-1}(\mathbf{f}) = \dim Z \cap X_{\mathbf{f}}  \geq \dim Z \cap X_{\mathbf{f}'} = \dim (\pi\circ i)^{-1}(\mathbf{f}').\end{equation}
 We restrict now to the subset $\pi^{-1}(U_{flat})$. Since $\pi: \pi^{-1}(U_{flat}) \rightarrow U_{flat}$ is flat, then $x \in \pi^{-1}(U_{flat}) \cap Z$ if and only if $X_{\pi(x)}$ is singular at $x$. In particular, if $\mathbf{f} \in U_{flat}$, then $$(\pi\circ i)^{-1}(\mathbf{f})=X_\mathbf{f} \cap Z= \mrm{Sing}(X_\mathbf{f}).$$ Recall now that if $O_\mathbf{f} \subset \overline{O_{\mathbf{f}'}}$ and $\mathbf{f} \in U_{flat}$ then $\mathbf{f}' \in U_{flat}$. As a consequence of Equation~\ref{eq:dimsing} we have, $$ \dim \mrm{Sing}(X_{\mathbf{f}})= \dim X_\mathbf{f} \cap Z  \geq \dim X_{\mathbf{f}'}\cap Z = \dim \mrm{Sing}( X_{\mathbf{f}'}).$$\endproof
\begin{rmk} Let $O_{\mrm{id}}$ denote the $G$-orbit of $(\mrm{id}, \dots, \mrm{id})$. Because of Theorem~\ref{thm:smooth}, a flat irreducible linear degeneration $\mathcal{F}l_\mathbf{d}^\mathbf{f} (\C^m)$ is singular if and only if ${\bf f} \in U_{flat, irr} \setminus O_{\mrm{id}}$.
\end{rmk}
Restricting the choices of possible dimension vectors, some additional informations about the singular locus can be deduced.
 \begin{theorem}\label{thm:singcodim3}Suppose that $d_{i+1}-d_i=1$ for all~$i<n$. If $X$ is a singular flat irreducible linear degeneration, then $\mrm{Sing}(X)$ is a closed subvariety of codimension 3. 
\end{theorem}
\proof 
 Recall that, because of Theorem~\ref{thm:reg2}, degenerations in $O_{\mathbf{r}^1}$ are regular in codimension 2. Consider now a degeneration $X=X_{\bf f}$ such that ${\bf f} \in U_{flat, irr}$. Since $U_{flat, irr}$ is the locus of orbits degenerating to~$O_{\mathbf{r}^1}$, as a consequence of Proposition~\ref{prop:closuredim} we have that $ \mrm{codim} Sing (X_{\bf f}) \geq 3.$ Because of Remark~\ref{rmk:upperdeg}, if ${\bf f} \in U_{flat, irr} \setminus O_{\mrm{id}}$, there exists $h$ such that $X_{\bf f}$  is a degeneration of $\mrm{Gr}_{\mathbf{d}}(M^h)$. Then by Proposition~\ref{prop:closuredim}, we have that 
$3 = \mrm{codim} \mrm{Sing}(\mrm{Gr}_{\mathbf{d}}(M^h)) \geq \mrm{codim} \mrm{Sing}(X )$
and the theorem is proved. \endproof
\begin{rmk} Observe that the previous theorem cannot be extended to all degenerations $\mathcal{F}l_\mathbf{d}^\mathbf{f} (\C^m)$ such that ${\bf f} \in U_{flat} \setminus O_{\mrm{id}}$. In fact, it is proved in \cite[Theorem 13]{CIFFFR17} that there exist flat non irreducible degenerations having singularities in codimension 1. 
\end{rmk}
According to Corollary~\ref{cor:codimsing}, the above theorem does not hold in general if the condition $d_{i+1}-d_i=1$ is not satisfied for all~$i<n$. Additionally, we provide an example showing that strict inequalities between the dimensions of singular loci can arise when the involved maps belong to different orbits.
 \begin{example}\label{ex:possiblecodimension}
For the partial flag variety $\mathcal{F}l_\mathbf{d} (\C^6)$ with $\mathbf{d}=(1,4)$, there are four isomorphism classes of flat irreducible linear degenerations. When identified with quiver Grassmannians, these classes correspond to representations
\[ M_{f}:= \C^6 \stackrel{f} \longrightarrow \C^6\]
where the rank of $f$ is $3, 4, 5,$ or $6$.
The case where $\rk f = 6$ corresponds to the classical partial flag variety itself, which is a smooth variety.

For the other cases, we can describe the dimension of their singular loci:
\begin{itemize}
    \item If $\rk f = 5$, Corollary~\ref{cor:codimsing} shows that the singular locus is a closed subvariety of dimension $4$.
    \item If $\rk f = 4$, the quiver Grassmannian $\mrm{Gr}_\mathbf{d}(M_f)$ is an irreducible Schubert quiver Grassmannian, which can be realized as a Schubert variety in a suitable partial flag variety. A direct computation reveals that the singular locus of this Schubert variety has dimension $6$.
    \item Analogously, if $\rk f = 3$, Theorem~\ref{thm:reg2} allows for a direct verification that the singular locus is a closed subvariety of dimension $8$.
\end{itemize}
\end{example}

We computed several additional examples that suggest a number of interesting questions, which we hope to explore in future work.
\begin{question}\label{q:1} Let $X$ be a well behaved quiver Grassmannian (for any Dynkin quiver). Is it true that $\mrm{Sing}(X)$ is a subvariety of codimension 3? \end{question}
\begin{question}\label{q:2} Is it true that if $X$ is a singular flat irreducible degeneration, then $\mrm{Sing}(X)$ is a variety of odd codimension? \end{question}
To provide context for the third question, let us fix an $n-1$ tuple $\mathbf{f}$ of endomorphisms and a dimension vector $\mathbf{d}$. Define
\[J_{\bf f} := \{ i \in \{1, \dots n-1\} \; | \; \rk f_i < m \}, \qquad D_{\bf f}=\min_{h \in J}\{d_{h+1}-d_h\}.\]
As a consequence of Corollary~\ref{cor:codimsing} and of Proposition~\ref{prop:closuredim}, the codimension of the singular locus of a flat irreducible degeneration  $\mathcal{F}l_\mathbf{d}^{\mathbf{f}} (\C^m)$ is an integer in $\{ 3, \dots, 2D_{\bf f}+1\}$. Set $D=\min_{h}\{d_{h+1}-d_h\}$.
\begin{question}\label{q:3} Consider the map $ \eta^c : U_{flat, irr} \rightarrow \mathbb{N} $ which assigns to each $n-1$-tuple $\mathbf{f}$ the codimension of the singular locus of $\mathcal{F}l_\mathbf{d}^{\mathbf{f}} (\C^m)$. Is it possible to describe the image of $\eta^c$? Does the image $  \eta^c( U_{flat, irr})$ always coincide with the set of odd integers in $ \{ 3, \dots 2D+1 \}$?  \end{question}

\subsection{Regularity of Irreducible Degenerations.}
Consider the generic irreducible degeneration $X=\mathcal{F}l_\mathbf{d}^\mathbf{f} (\C^m)$.
As a result of Corollary~\ref{cor:euler}, we have $X=X_1 \times \dots \times X_k$ where each $X_i$ is a flat irreducible linear degeneration.
The results from Section~\ref{subsec:regularity} enable us to deduce some general results about the regularity of irreducible degenerations.
\begin{theorem}\label{thm:irreg2} Every irreducible degeneration is regular in codimension 2.
\end{theorem}
\proof First, observe that the singular locus of $X=X_1 \times \dots \times X_k$ is given by \[\mrm{Sing}(X)= \bigcup_{i=1}^k  \Sigma_i, \]
where \[ \Sigma_i  = X_1 \times \dots \times X_{i-1} \times \mrm{Sing}(X_i) \times X_{i+1} \times \dots X_k. \]
To see this, consider $p=(p_1, \dots, p_k) \in X$, where $p_i \in X_i$. Since $$\dim T_p(X) = \dim T_{p_1}(X_1)+ \dots + \dim T_{p_k}(X_k),$$ a point $p$ is singular in $X$ if and only if at least one $p_i$ is a singular point in $X_i$.
Consequently, \begin{equation}\label{eq:codimensionproduct}\mrm{codim}_X\mrm{Sing}(X)=\min_i \mrm{codim}_{X_i} \mrm{Sing} (X_i).\end{equation}According to Theorem~\ref{thm:reg2}, every $X_i$ is regular in codimension 2. Thus, Formula~\ref{eq:codimensionproduct} implies the claim.\endproof 
If the dimension vector satisfies $d_{i+1}-d_i=1$ for all $i \in \{1, \dots, n\}$, more precise dimension estimates for the singular locus hold.
\begin{theorem}Suppose that $d_{i+1}-d_i=1$ for all $i \in \{1, \dots, n\}$. If $X$ is a singular irreducible linear degeneration, then $\mrm{Sing}(X)$ is a subvariety of $X$ of codimension 3. 
\end{theorem}
\proof As observed in the proof of Theorem~\ref{thm:irreg2}, $X=X_1 \times \dots \times X_k$ is singular if and only if at least one of the $X_i$ is singular. Due to Theorem~\ref{thm:singcodim3}, if $X_i$ is a singular flat irreducible degeneration, then $\mrm{Sing}(X_i)$ is a closed subvariety of codimension 3.
The assertion now follows from Formula~(\ref{eq:codimensionproduct}).
\endproof

\end{document}